\documentclass{amsart}
\input xypic
\usepackage{amssymb}   % For Latex2e
\usepackage{amsmath}
\usepackage{amsthm}
\usepackage{times}
\newtheorem{Theorem}{Theorem}
\newtheorem{Lemma}[Theorem]{Lemma}
\newtheorem{Proposition}[Theorem]{Proposition}
\newtheorem{Corollary}[Theorem]{Corollary}

\theoremstyle{definition}

\theoremstyle{definition}
\newtheorem{Definition}[Theorem]{Definition}
\newtheorem{Remark}[Theorem]{Remark}
\newtheorem{Example}[Theorem]{Example}

\numberwithin{Theorem}{section} 
\numberwithin{equation}{section}

\newcommand{\Pic}{\operatorname{Pic}}
\newcommand{\Bl}{\operatorname{Bl}}

\newcommand{\Tx}{\operatorname{T}_x}

\renewcommand{\O}{{\mathcal O}}

\newcommand{\I}{{\mathcal I}}

\newcommand{\Proj}{{\mathbb P}}

\newcommand{\p}{{\mathbb P}}

\newcommand{\codim}{\operatorname{codim}}

\renewcommand{\H}{\operatorname{H^0}}
\newcommand{\Hi}{\operatorname{H}}

\newcommand{\map}{\dasharrow}

\newcommand{\Num}{\operatorname{Num}}
\newcommand{\NS}{\operatorname{NS}}

\newcommand{\defect}{\operatorname{def}}

\newcommand{\ho}{\operatorname{h^0}}
\def\S{Section~}

\begin{document}
\title[Quadratic entry locus varieties ]{Varieties with quadratic entry locus, $ I$}
\author{Francesco Russo}
\address{Dipartimento di Matematica e Informatica, Universit\` a degli Studi di Catania,
Viale A. Doria 6, 95125 Catania, Italia}
\email{frusso@dmi.unict.it}
\thanks{Partially supported by CNPq
(Centro Nacional de Pesquisa), grant 308745/2006-0, and by
PRONEX/FAPERJ--Algebra Comutativa e Geometria Algebrica}

\maketitle

%\tableofcontents
\section*{Introduction}
We shall introduce and study  quadratic entry locus varieties, a
 class of projective algebraic varieties whose extrinsic
and intrinsic geometry is very rich.
\thispagestyle{empty}

Let us recall that, for   an irreducible non-degenerate variety
$X\subset\p^N$ of dimension $n\geq 1$,  {\it the secant defect of
$X$}, denoted by $\delta(X)$, is
 the difference between the
expected dimension and the effective dimension of the secant variety
$SX\subseteq\p^N$ of $X$, that is $\delta(X)=2n+1-\dim(SX)$. This is
an important projective invariant measuring the dimension of the
{\it entry locus} $\Sigma_p\subseteq X$ described by the  points on
$X$ spanning secant lines passing through a general point $p\in SX$;
see \S 1.

Many  examples appearing in different settings
suggested the definition of quadratic entry locus manifold of type
$\delta$, briefly $QEL$-manifold of type $\delta$. These are smooth
varieties $X\subset\p^N$ for which $\Sigma_p$ is a smooth quadric
hypersurface of dimension $\delta=\delta(X)$, whose linear span in
$\p^N$ is the locus of secant lines to $X$ passing through $p$,
$p\in SX$ general. We also consider $LQEL$-manifolds of type $\delta$, that is smooth varieties $X\subset\p^N$ for which $\Sigma_p$
is the union of smooth quadric hypersurfaces of dimension $\delta=\delta(X)$, $p\in SX$ general; see \S 1 for more details. 
When $\delta=0$ and $N> 2n+1$ the above
conditions  do not impose particular geometric  restrictions on $X$. On the
contrary $QEL$-manifolds  $X\subset\p^{2n+1}$ of type $\delta=0$
are linearly normal and rational, see \cite{CMR}, while in \cite{IR2} it is proved
that every $QEL$-manifold of type $\delta>0$ is rational. The notion of
$(L)QEL$-manifold  was also motivated by the  remark  that  a
lot of secant defective smooth varieties with special geometric properties
and/or with  extremal tangential behaviour are $(L)QEL$-manifolds:   varieties defined by quadratic equations
having enough linear syzygies, for example satisfying condition
$N_2$ of Green (Proposition \ref{K2}); homogeneous varieties, secant
defective or not; Scorza Varieties (and in particular Severi
Varieties; this property being an essential ingredient for their
classification, cf.\ \cite{Zak1} and \S 3); centers of special Cremona transformations of type
$(2,d)$ (\cite{ESB} and \S \ref{Cremonasection}), varieties whose dual variety
is small (\cite{Ein} and \cite{IR2}) and which are not hypersurfaces. Furthermore, if
$n=2,3,$ any smooth secant defective variety $X\subset\p^N$ with
$SX\subsetneq\p^N$ is a $LQEL$-manifold; see \cite{Sev,Scorza1,Fujitasec}.

By definition  $LQEL$-manifolds of type $\delta>0$
are   very
special examples of rationally connected varieties \cite{KMM,Kollar,Debarre}. Gaetano Scorza was the first who realized
the link between secant defective varieties and rational
connectedness in the pioneering papers \cite{Scorza1,Scorza2}, where, {\it ante litteram}, the condition of
rational connectedness by conics appears for {\it le variet\' a
della $(n-1)$-esima specie, o dell' ultima specie} in \cite[pp.\
252--253]{Scorza2}. These definitions are explained in detail in
\cite[Section 2.3]{Russo3} and also inspired some  results  in   \cite{IR}  and in  \cite{IR2}.

 Here we develop the theory of $LQEL$-manifolds of
type $\delta\geq 1$ by studying the geometry of the family of conics
and, for $\delta\geq 2$, of lines produced by the quadratic entry
loci via the modern tools of deformation theory of rational curves
on a manifold and via their parameter spaces. For $\delta\geq 2$,
the lines passing through a general point $x\in X$  describe a
smooth, not necessarily irreducible, variety
$Y_x\subset\p((\mathbf{T}_{x,X})^*)=\p^{n-1}$, defined in
Proposition \ref{Yx}.

The most important results for $LQEL$-manifolds of type $\delta\geq
2$ are consequences of the study of the projective geometry of the
subvariety $Y_x\subset\p^{n-1}$, especially for $\delta\geq 3$ when
$Y_x$ is also irreducible. 
This is not surprising since the family of
lines on such an $X\subset\p^N$ is the {\it minimal covering family
of rational curves} in the sense of Mori; see \cite{Mori1}, \cite{HM0}, \cite{HM1}, \cite{HM2}, the surveys 
\cite{Hwang}, \cite{Hwang2} and the recent \cite{HK} for many interesting connections between the projective geometry of analogous of $Y_x$, the so called
{\it variety of minimal rational tangents}, 
and geometrical properties of Fano manifolds.

With regard to the projective geometry of $Y_x\subset\p^{n-1}$, we
prove a key result asserting that, for $LQEL$-manifolds of type
$\delta\geq 3$, the variety $Y_x\subset\p^{n-1}$ is a $QEL$-manifold
of type $\delta-2$ such that $SY_x=\p^{n-1}$, see Theorem
\ref{quadrichetangenti}, part (d). This allows the inductive definition of
some varieties naturally attached to a $LQEL$-manifold of type
$\delta\geq 3$, see \S \ref{qualitative}, which can be considered as
 algebro-geometric analogues of the successive projective
differential forms of an arbitrary irreducible projective variety;
\cite{GH,IL}. If $r_X=\sup\{r\in \mathbb{N} :
\delta\geq 2r+1\}$,
 for every $k=1, \ldots, r_X-1$, we set
$X^k=Y_x(X^{k-1}),$ where $x\in X^{k-1}$ is a general point and
where $X^0=X$ (see \S \ref{qualitative} for
more precise definitions). Thus for every $k=1,\ldots, r_X-1$, the
variety $X^k$ is a $QEL$-manifold of type $\delta(X^k)=\delta-2k\geq
3$. Performing the calculations on the dimension of the $X^k$'s we
arrive at the surprising {\it Divisibility Property} that $2^{r_X}$
divides $n-\delta$ and at other strong restrictions,  Theorem~\ref{fundamental}. 
This result can be considered as a generalized secant
analogue of the famous  Landman Parity Theorem for the dual
defect $\defect(X)$ of a variety $X\subset\p^N$: $n\equiv\defect(X)$
(mod $2$), if $\defect(X)>0$; see \cite[Theorem 2.4]{Ein}. The repeated
applications of the Divisibility Theorem in the paper will allow
unitary and simple proofs of  many important theorems regarding projective varieties
with very special geometric  properties and it was also applied, together with its consequences proved here,  in 
\cite[Section 4]{Fu} to obtain a new Linear Normality Bound for $LQEL$-manifolds  and a new characterization of smooth quadric hypersurfaces. 

For example we can easily classify $LQEL$-manifolds of type $\delta>\frac{n}{2}$ in Corollary \ref{deltatranensudue},
 obtaining only five known
examples and answering a problem posed in \cite[0.12.6]{KS}. We  also provide a two line proof 
that $LQEL$-manifolds of type $\delta=\frac{n}{2}$ appear only in
dimension 2, 4, 8 or 16, see the proof of Corollary
\ref{nminore16}. In particular this last result  implies that Severi varieties appear
only in these dimensions. This was the key and most difficult point for getting their classification, 
which after this step is quite easy at least today, see the almost self-contained and
short proof of the classification in Corollary \ref{nminore16}. The classification of Severi varieties was the main
result  of \cite{Zak3}  (see also \cite{LVdV,  Zak1, Lan2, Chaput}) and  produced a great impact for its deep connections with Hartshorne Conjecture
and also
with other areas of mathematic such as the classification of composition algebras over a field.
Here, via our approach, we generalize this to  the classification of $LQEL$-manifolds of type $\delta=n/2$, see {\it loc. cit.},
and we present an interesting connection between the classification of Severi Varieties and Hartshorne Conjecture which,
as far as we know, was not noticed before, see Remark \ref{HartYx}.

In \S \ref{Cremonasection}, using also \cite{IR},  we obtain some strong
restrictions for the existence of special Cremona transformations
$\phi:\p^N\map\p^N$ of type $(2,d)$, whose  base loci  are linearly
normal $QEL$-manifolds of type $\delta=\frac{n+2}{d}-1$ with  notable 
 properties, as shown in Proposition
\ref{adjunctionCremona}.  
 Our  applications to special Cremona transformations
concern the complete classification of
those  of type $(2,3)$, Corollary~\ref{special23}, 
of those of type $(2,5)$, Corollary~\ref{special25},
 and
of all special Cremona transformation in $\p^{2n+1}$, Corollary \ref{2npiu1},
where $n$ is the dimension of the base locus.

\section{Definitions, preliminary results and
examples}\label{prel} We  work over the field of complex
numbers, unless otherwise stated. 
 Let $X\subset\p^N$ be a closed
irreducible subvariety, which we will always suppose to be
non-degenerate. From now on $n=\dim(X)\geq 1$ and $N\geq n+1$.

Let $X\subset \Proj^N$ be an irreducible, non-degenerate projective
variety of dimension $n$. Let
\[SX=\overline{\textstyle{\bigcup\limits _{{\scriptstyle x\neq y }\atop {\scriptstyle x,y\in
X}}}\langle x,y\rangle }\subseteq\p^N\] be the {\it secant variety to $X$}.

Clearly $\dim(SX)\leq\min\{N,2n+1\}$. If the equality
$\dim(SX)=2n+1$ holds, then $X\subset\p^N$ is said to be {\it
non-defective}. If $\dim(SX)<2n+1$, then $X\subset\p^N$ is said to
be {\it secant defective}, or simply {\it defective} and
$\delta(X):=2n+1-\dim(SX)$ is called the {\it secant defect of
$X\subset\p^N$}.

For $p\in SX\setminus X$, the closure of the locus of couples of
distinct  points on $X$ spanning secant lines passing through $p$ is
called {\it the entry locus of $X$ with respect to $p\in SX$} and it
will be indicated by $\Sigma_p(X)$. The closure of the secant lines
to $X$ passing through $p$ is a cone over $\Sigma_p(X)$, let us call
it $C_p(X)$.  If $X\subset\p^N$ is smooth, then
$\Sigma_p(X)=C_p(X)\cap X$ as schemes for $p\in SX$ general, see for
example \cite[Lemma 4.5]{FR}. Moreover it is easy to see that for
$p\in SX$ general,  $\Sigma_p(X)$ is  equidimensional of dimension
equal to $\delta(X)$. Thus $\dim(C_p(X))=\delta(X)+1$. In general
$\Sigma_p(X)$ (and hence a fortiori $C_p(X)$) is not irreducible.

For an irreducible variety $X\subset\p^N$ and for $x\in X$, the
(embedded) projective tangent space at $x\in X$ is denoted by $\Tx
X$ and the affine (or Zariski) tangent space at $x\in X$ by~$\mathbf{T}_xX$.

\begin{Definition} (cf. also \cite{Zak1, KS, IR}) An irreducible, non-degenerate projective variety
$X\subset\p^N$ is said to be a {\it local quadratic entry locus
variety of type $\delta\geq 0$}, briefly an {\it LQEL-variety of type
$\delta$}, if, for general $x,y\in X$ distinct points and
 for general $p\in \langle x,y\rangle \subseteq SX$, the union of the irreducible components of
the entry locus of $p$ passing through $x$ and through $y$ is  a quadric
hypersurface of dimension $\delta=\delta(X)$ in the given embedding
$X\subset\p^N$. Equivalently an  $LQEL$-variety of type $\delta\geq 0$ is an irreducible projective variety $X\subset\p^N$
if through two general points
there passes a  quadric hypersuface of dimension $\delta=\delta(X)$ contained in $X$.

An irreducible projective variety $X\subset\p^N$ is said to be a {\it quadratic entry locus
variety of type $\delta\geq 0$}, briefly a {\it $QEL$-variety of
type $\delta$}, if for general $p\in SX$ the entry locus
$\Sigma_p(X)$ is a quadric hypersurface of dimension
$\delta=\delta(X)$.

An irreducible, non-degenerate projective variety $X\subset\p^N$ is
said to be a {\it conic-con\-nected variety}, briefly a {\it
$CC$-variety}, if through two general points of $X$ there passes an
irreducible conic contained in $X$. 

If $X\subset\p^N$ is also
smooth, we shall use the terms $LQEL$-manifold, $QEL$-manifold,
respectively  {\it CC}-manifold.
\end{Definition}
 
\begin{Lemma}\label{1.4} Let $X$ be an LQEL-manifold with  $\delta(X)>0$ and let $x,y \in X$ be general points. There is a unique quadric hypersurface of dimension $\delta$, say $Q_{x,y}$, passing through $x,y$ and contained in $X$. Moreover, $Q_{x,y}$ is irreducible.  
\end{Lemma}

\begin{proof} 
Uniqueness follows from the fact that the general entry locus passing through two general points is smooth at these points by Terracini Lemma (see e.g. \cite[Proposition 3.3]{HR}). To see that $Q_{x,y}$ is irreducible, we may assume $\delta=1$ by passing to general hyperplane sections (see Proposition
 \ref{deltan} below). 
Assume first that $X$ is covered by lines passing through $x$. Being also smooth, $X$ is a linear space,
so it is not an LQEL-manifold. 
Otherwise, after suitable normalization, the family of conics through $x$ is generically smooth and the result follows.  
\end{proof}

A monodromy argument shows that for an
$LQEL$-variety the general entry locus $\Sigma_p(X)$ is a union of
quadric hypersurfaces of dimension $\delta(X)$. Moreover, the general entry locus
of a  $QEL$-manifold, or each irreducible component of
the general  entry locus of an $LQEL$-manifold, is a smooth quadric hypersurface, 
see the arguments in
 \cite[pp.\ 964--966]{FR}.
\label{smoothnessEL}

A $QEL$-variety of type $\delta$ is clearly an $LQEL$-variety of type
$\delta$ and an $LQEL$-variety of type $\delta\geq 1$ is a
$CC$-variety but the converses are not true. Every variety $X\subset\p^N$ with $\delta(X)=0$ is a
$LQEL$-variety of type $\delta=0$, while $QEL$-varieties of type
$\delta=0$ are those for which through a general point of $SX$ there
passes a unique secant line. This last condition is especially
relevant for $N=2n+1$; see \cite{CMR} for example. One
constructs examples of $LQEL$-manifolds which are not
$QEL$-manifolds by projecting isomorphically a $QEL$-manifold
$X\subset\p^N$ of type $\delta\geq 1$ and with  $N>2n+1-\delta$ into~$\p^{2n+1-\delta}$.
Indeed, consider a general $L=\p^{N-2n+\delta-2}$ in $\p^N$ and let $\pi_L:\p^N\map\p^{2n+1-\delta}$ be the projection from $L$. Then $SX\cap L=\emptyset$ and the
the restriction of $\pi_L$ to  $SX$ is a finite surjective morphism onto  $\p^{2n+1-\delta}$ of degree $d\geq 2$ equal to $\deg(SX)$, while
the restriction of  $\pi_L$ to $X$ is an isomorphism onto $\overline X\subset\p^{2n+1-\delta}$. Thus a general $q\in S\overline X=\p^{2n+1-\delta}$
has exactly $d\geq 2$ preimages via $\pi_L$ on $SX$, let us say $p_1,\ldots,p_d$, which are also general points on $SX$. The entry locus $\Sigma_q(\overline X)$
contains as irreducible components of dimension $\delta=\delta(\overline X)$ the projection of the $d$ quadric hypersurfaces $\Sigma_{p_i}(X)$. Moreover, it is not difficult to see that the irreducible components of $\Sigma_q(\overline X)$
of dimension $\delta=\delta(\overline X)$ are exactly the quadric hypersurfaces $\pi_L(\Sigma_{p_i}(X))$, $i=1,\ldots, d$.

Let us collect some  consequences of the definitions in the
following Proposition, whose  proof is left to the reader, see also
\cite[Proposition 1.3]{IR2}.

\begin{Proposition}\label{deltan}
Let $X\subset\p^N$ be an irreducible, non degenerate 
projective variety and let, if it exists, $X'\subset\p^M$, $M\leq
N-1$, be an isomorphic projection of $X\subset\p^N$. Then
\begin{enumerate}
\item[i)] If $SX=\p^N$ and if $X\subset\p^N$ is a $QEL$-manifold,
then $X\subset\p^N$ is linearly normal;
\item[ii)] $X'\subset\p^M$ is an $LQEL$-manifold if and only if
$X\subset\p^N$ is an $LQEL$-manifold;
\item[iii)] If $X\subset\p^N$ is a(n) $(L)QEL$-manifold of type $\delta\geq
1$, then  a general hyperplane section
$\widetilde{X}\subset\p^{N-1}$ is a(n) $(L)QEL$-manifold of type
$\widetilde{\delta}=\delta-1$.
\item[iv)] $X\subset\p^N$ is an $LQEL$-variety of type $\delta=n$ if and only if $N=n+1$ and $X\subset\p^{n+1}$ is
a  quadric hypersurface.
\end{enumerate}
\end{Proposition}

We describe a lot of examples of $QEL$-manifolds to see their
ubiquity among projective varieties with very special geometric
properties, as cited in the introduction. Another interesting class,
not contained in the one described below, is considered in \S
\ref{Cremonasection}.

\begin{Proposition}\label{K2}{\rm (\cite{Ve}, \cite{HKS})}
 A smooth non-degenerate variety $X \subset \p^N$, scheme theoretically defined by quadratic equations whose Koszul syzygies are generated by linear ones, is a QEL-manifold.
\end{Proposition}

 There are several possible equivalent
definitions of the projective second fundamental form
$|II_{x,X}|\subseteq\p(S^2(\mathbf{T}_xX))$ of an irreducible
projective variety $X\subset\p^N$ at a general point $x\in X$, see
for example \cite[3.2 and end of Section~3.5]{IL}. We shall use the
one related to tangential projections, as in \cite[remark
3.2.11]{IL}.

 Suppose $X\subset\p^N$ is non-degenerate, as always,
 let $x\in X$ be a general point and
consider the projection from $\Tx X$ onto a disjoint $\p^{N-n-1}$:
\begin{equation}\label{tangentdef}
\pi_x:X\map W_x\subseteq\p^{N-n-1}.
\end{equation}
The map $\pi_x$ is associated to the linear system of hyperplane
sections cut out by hyperplanes containing $\Tx X$, or equivalently
by the hyperplane sections of $X\subset\p^N$ singular at $x$.

Let $\phi:\Bl_xX\to X$ be the blow-up of $X$ at $x$, let
\[E=\p((\mathbf{T}_xX)^*)=\p^{n-1}\subset\Bl_xX\] be the exceptional
divisor and let $H$ be a hyperplane section of $X\subset\p^N$. The
induced rational map $\widetilde{\pi}_x:\Bl_xX\map\p^{N-n-1}$ is
defined as a rational map along $E$ since $X\subset\p^N$ is not a linear space; see
for example the argument in \cite[2.1 (a)]{Ein}. The restriction of
$\widetilde{\pi}_x$ to $E$ is given by a linear system in
$|\phi^*(H)-2E|_{|E}\subseteq|-2E_{|E}|=|\O_{\p((\mathbf{T}_xX)^*)}(2)|=\p(S^2(\mathbf{T}_xX))$.

\begin{Definition}
 The {\it second fundamental form
$|II_{x,X}|\subseteq\p(S^2(\mathbf{T}_xX))$ of an irreducible
non-degenerate variety $X\subset\p^N$ of dimension $n\geq 2$ at a
general point $x\in X$} is the non-empty linear system of quadric
hypersurfaces in $\p((\mathbf{T}_xX)^*)$ defining the restriction of
$\widetilde{\pi}_x$ to $E$.
\end{Definition}

Clearly  $\dim(|II_{x,X}|)\leq N-n-1$ and
$\widetilde{\pi}_x(E)\subseteq W_x\subseteq\p^{N-n-1}$. From this
point of view the base locus on $E$ of the second fundamental form
$|II_{x,X}|$ consists of {\it asymptotic directions}, i.e. of
directions associated to lines having a contact of order at least two with
$X$ at $x$. For example when $X\subset\p^N$ is defined by equations of
degree at most two,  the base locus of the second fundamental form
consists of points giving tangent lines contained in $X$ and passing
through $x$ so that it is exactly the locus of lines through $x$ and contained in $X$,
$Y_x\subset\p((\mathbf{T}_xX)^*)=E$, which will be defined in
Proposition \ref{Yx}.

\begin{Lemma}\label{2.1}
Let $X\subset \p^N$ be a smooth irreducible non-degenerate variety, and assume $\delta(X)>0$. The irreducible components of the closure of a general fibre of $\pi_x$ are not~linear.  
\end{Lemma}
\begin{proof} 
We may assume by passing to linear sections $\delta(X)=1$. Let $l$ be a line, passing through a general point $y\in X$, which is an irreducible component of the closure of a general fibre of $\pi_x$. 
By Terracini Lemma $T_x(X)\cap T_y(X)$ is a point, say $p_{x,y}$.
Since $l\subset \langle T_{x}(X), y\rangle\cap T_y(X)$, $p_{x,y} \in l$. 
By symmetry there is also a line $l'$ in  
$ \langle T_{y}(X), x\rangle \cap X$, $x\in l'$, $p_{x,y}\in l'$. So, $l\cup l'$ is a conic contained in the plane $\langle x,y,p_{x,y}\rangle $ passing through $x,y$. Reasoning as in the proof of Lemma \ref{1.4} we find a contradiction. 
\end{proof}

\section{Qualitative  properties of
$LQEL$-manifolds}\label{qualitative} We describe the conics
naturally appearing  on $LQEL$-manifolds of type $\delta>0$ and
relate them to intrinsic invariants; see also  \cite{IR2}.

\begin{Theorem}\label{conics} Let $X\subset\p^N$ be an
 $LQEL$-manifold of type $\delta\geq 1$. Then:
\begin{enumerate}

\item The variety $X$ is  rationally connected, so that it is a
simply connected manifold such that  $
\H(\Omega^{\otimes m}_X)=0$ for every
$m\geq 1$ and $\Hi^i(\O_X)=0$ for every $i>0$.

\item There exists on $X$ an irreducible family of  conics
$\mathcal{C}$ of dimension $2n+\delta-3$, whose general member is
smooth. This family describes an open subset of an irreducible
component of the Hilbert scheme of conics on $X$.
\item Given a general point $x\in X$, let $\mathcal{C}_x$ be the  family of  conics in
$\mathcal{C}$ passing through $x$. Then $\mathcal{C}_x$ has
dimension $n+\delta-2$, equal to the dimension of the irreducible
components of $\mathcal{C}_x$ describing dense subsets of $X$.

\item Given two general points $x,y\in X$, the locus  $Q_{x,y}$ of the family $\mathcal{C}_{x,y}$
of smooth conics in $\mathcal{C}$ passing through $x$ and $y$ is
a smooth quadric hypersurface of dimension
$\delta$. The family $\mathcal{C}_{x,y}$ is irreducible and of
dimension $\delta-1$. 

\item For a general  conic $C\in \mathcal{C}$,
$-K_X\cdot C=n+\delta.$

\item A general conic $C\in \mathcal{C}_x$ intersects $\Tx X$ only
at $x$. Moreover the tangent lines to smooth conics contained in $X$
and passing through $x\in X$ describe an open  subset of
$\p((\mathbf{T}_xX)^*)$.
\end{enumerate}
\end{Theorem}
\begin{proof} The variety $X$ is clearly rationally connected. The
conclusions of part (1) are contained in \cite[Proposition 2.5]{KMM} and
also in \cite[Corollary
4.18]{Debarre}.

Part (4) is the definition of an $LQEL$-variety. Indeed, the plane
spanned by every conic through $x$ and $y$ contains the line
$\langle x,y\rangle $, which is a general secant line to $X$. Thus a conic through
$x$ and $y$  is  contained in the entry locus of every $p\in \langle x,y\rangle $
not on $X$, so that for $p\in \langle x,y\rangle $ general it is contained in the
smooth quadric hypersurface $Q_{x,y}$, the unique irreducible component of
$\Sigma_p$ passing through $x$ and $y$.

Let us prove parts (3) and  (5). Fixing two general points $x,y \in
X$ there exists a smooth quadric hypersurface $Q_{x,y}\subset
X\subset\p^N$ of dimension $\delta\geq 1$ through $x$ and $y$. Thus
there exists a smooth conic passing through $x$ and a general point
$y\in X$. In particular, there exists an irreducible family of smooth
conics passing through $x$, let us say $\mathcal{C}^1_x$, whose
members describe a dense subset of $X$. Since the number of
irreducible components of the Hilbert scheme of conics in $X$ is
finite, we get that
\[T_{X|C}\simeq{\textstyle \bigoplus\limits_{i=1}^n}\O_{\p^1}(a_i)\]
is ample for $C\in \mathcal{C}^1_x$ general; see for example
\cite[II.3.10.1]{Kollar}. Hence  $a_i>0$ for every $i=1,\ldots, n$,
and $N_{C/X}$ is ample being a quotient of $TX_{|C}$. Thus
$\mathcal{C}^1_x$ is smooth at $C\in \mathcal{C}^1_x$ and of
dimension
\[H^0(N_{C/X}(-1))=H^0(T_{X|C}(-1))-2=-2+\sum_{i=1}^na_i=-K_X\cdot C-2.\] Let
$\widetilde{\mathcal{C}_x}$ be the universal family over
$\mathcal{C}_x$ and let $\pi:\widetilde{\mathcal{C}_x}\to X$ be the
tautological morphism. By part (4) and the theorem on the dimension
of the fibers we get that
$\dim(\widetilde{\mathcal{C}^1_x})=n+\delta-1$. Thus
$$n+\delta-2=\dim(\mathcal{C}^1_x)=-K_X\cdot C-2$$
and part (3) and (5) are proved. Part (2) now easily follows.

Let $C\subset X$ be a general smooth conic passing through $x\in X$.
In particular $C$ is not completely contained in $\Tx X\cap X$ since
the conics through $x$ cover $X$. Suppose there exists a point $z\in
\Tx X\cap C$, $z\neq x$. The line $\langle x,z\rangle $ is contained in $\Tx X$
and also in the plane generated by $C$, let us say $\Pi=\langle C\rangle $. Since
the conic $C$ is smooth at $x$ and since $C$ is not contained in
$\Tx X$, we deduce $\langle x,z\rangle \subseteq \Tx X\cap\Pi=T_xC$, that is
$\langle x,z\rangle =T_xC$. Then the line $\langle x,z\rangle $ would cut $C$ in at least $3$
points, counted with multiplicities, which is clearly impossible.

 We saw that $TX_{|C}\simeq\bigoplus\limits_{i=1}^n\O_{\p^1}(a_i)$, with
$a_i>0$ for $C\in \mathcal{C}_x$ general. Consider the map
$\tau_x:\mathcal{C}_x\map\p((\mathbf{T}_xX)^*)$, which associates to
a conic in $\mathcal{C}_x$ its tangent line at $x$. The closure of
the image of $\tau_x$ in $\p((\mathbf{T}_xX)^*)$ has dimension
$n-1$. Indeed by \cite[II.3.4]{Kollar} the dimension of
$\tau_x(\mathcal{C}_x)$ at $\tau_x(C)$ is equal to $\#\{a_i>0\}-1$,
finishing the proof.
\end{proof}

On an $LQEL$-manifold of type $\delta\geq 2$ there are also lines
coming from the entry loci and we proceed to investigate them. 
The following result is
essentially well known, see also \cite[Proposition 1.5]{Hwang};
we recall its proof for reader's convenience.

\begin{Proposition}\label{Yx}
Let $X\subset\p^N$ be a smooth irreducible variety. Then:
\begin{enumerate}
\item The Hilbert scheme of lines passing through a general point
$x\in X$, if not empty, is smooth and can be identified with  a
smooth not necessarily irreducible variety
$Y_x\subseteq\p^{n-1}=\p((\mathbf{T}_xX)^*)$.
\item If $Y_x^j$, $j=1,\ldots,m$, $m\geq 2$, are the irreducible components
of $Y_x$, then  we have \[\dim(Y_x^l)+\dim(Y_x^p)\leq n-2 \quad \mbox{for every } l\neq
p.\] 
\end{enumerate}
\end{Proposition}
\begin{proof}We argue essentially as in the proof of
Theorem \ref{conics}. Let $L$ be a line through  the very general
point $x\in X$. Every such line $L$ is free, see for example
\cite[II.3.11]{Kollar}, so that
$TX_{|L}=\bigoplus\limits_{i=1}^n\O_{\p^1}(a_i(L))$, with $a_n(L)\geq
a_{n-1}(L)\geq\cdots\geq a_1(L)\geq 0$. Moreover $a_n(L)\geq 2$
because $T\p^1$ is a subbundle of $TX_{|L}$. On the other hand,
$T\p^N_{|L}=\O(2)\oplus\O(1)^{N-1}$ contains $TX_{|L}$ as a
subbundle so that $a_n(L)=2$ and $1\geq a_{n-1}(L)\geq\cdots\geq
a_1(L)\geq 0$, i.e.\ the arbitrary line $L$ is  a  {\it standard}  (or
{\it minimal}) curve in the sense of Mori Theory. It follows that
the map which associates to each line through $x$ its tangent
direction is a closed embedding so that we can identify the Hilbert
scheme of lines through $x$ with a variety
$Y_x\subset\p^{n-1}=\p((\mathbf{T}_xX)^*)$, which is smooth. Indeed,
$N_{L/X}=\bigoplus\limits_{j=1}^{n-1}\O_{\p^1}(b_j(L))$ with $b_j(L)\geq 0$
for every $j=1,\ldots, n-1$, being the quotient of a locally free
sheaf generated by global sections. Therefore $h^1(N_{L/X}(-1))=0$
and $Y_x$ is smooth at the point corresponding to $L$. Since $L$ was
an arbitrary line through $x$, $Y_x$ is smooth. The conditions on
the dimension of two irreducible components simply say that these
components  cannot intersect in $\p^{n-1}$.
\end{proof}

Now we prove a fundamental result on the geometry of lines on an
$LQEL$-manifold, which, via the study of the projective geometry of
$Y_x\subset\p^{n-1}$ and of its dimension, will yield 
significant obstructions for the existence of $LQEL$-manifolds of
type $\delta\geq 3$. The most relevant part for future applications
is part (4), d). Part (1) holds more generally for every smooth secant
defective variety.

\begin{Theorem}\label{quadrichetangenti}
Suppose that $X\subset\p^N$ is  an $LQEL$-manifold of type $\delta$.
\begin{enumerate}

\item {\rm (\cite[p.\ 282, Opere Complete, vol. I]{Scorza2})} If
$\delta\geq 1$, then $$\widetilde{\pi}_x:\p((\mathbf{T}_xX)^*)\map
W_x\subseteq\p^{N-n-1}$$ is dominant, so that
$\dim(|II_{x,X}|)=N-n-1$ and $N\leq\frac{n(n+3)}{2}$.

\item If $\delta\geq 2$, the smooth, not necessarily irreducible,
variety $Y_x\subset\p^{n-1}$ is non-degenerate and it consists of
irreducible components of the base locus scheme of $|II_{x,X}|$.
Moreover, the closure of the irreducible component of a general
fiber of $\widetilde{\pi}_x$  passing through a general point $p\in
\p((\mathbf{T}_xX)^*)$ is a linear space $\p^{\delta-1}_p$, cutting
scheme-theoretically $Y_x$ in a quadric hypersurface of dimension
$\delta-2$.

\item If $Y_x\subset\p((\mathbf{T}_xX)^*)$ is irreducible and if
$\delta\geq 2$, then $SY_x=\p((\mathbf{\mathbf{T}}_xX)^*)$ and
$Y_x\subset\p((\mathbf{T}_xX)^*)$ is a $QEL$-manifold of type
$\delta(X)-2$.

\item If $\delta\geq 3$, then

\begin{enumerate}
\item[{\rm a)}] $\Pic(X)\simeq\mathbb{Z}\langle \O_X(1) \rangle $.

\item[{\rm b)}]   For any line $L\subset X$,
$-K_X\cdot L=\frac{n+\delta}{2},$
so that
$i(X)=\frac{n+\delta}{2},$
where $i(X)$ is the index of the Fano manifold $X$.
 In particular
$n+\delta\equiv 0\;\;({\rm mod}~2),$ that is $n\equiv\delta\;\;({\rm mod}~2)$.

\item[{\rm c)}] There exists on $X$ an irreducible family of lines  of
dimension $\frac{3n+\delta}{2}-3$ such that for a general $L$ in
this family
\[TX_{|L}=\O_{\p^1}(2)\oplus\O_{\p^1}(1)^{\frac{n+\delta}{2}-2}\oplus\O_{\p^1}^{\frac{n-\delta}{2}+1}.\]

 \item[{\rm d)}] If $x\in X$ is general, then  $Y_x\subset\p((\mathbf{T}_xX)^*)$ is a 
  $QEL$-manifold of dimension $\frac{n+\delta}{2}-2$, of
type $\delta(X)-2$ and such that $SY_x=\p((\mathbf{T}_xX)^*)$;

\end{enumerate}
\end{enumerate}
\end{Theorem}
\begin{proof} Part (1) is classical and as we said above holds  for every 
smooth secant defective variety. Since its proof is self-contained and elementary
for $LQEL$-varieties, we include it for the reader's convenience. It
suffices to show that, via the restriction of $\widetilde{\pi}_x$,
the exceptional divisor $E=\p((\mathbf{T}_xX)^*)$ dominates
$W_x\subseteq\p^{N-n-1}$. Take a general point $y\in X$. By part (6) of
Theorem \ref{conics}, there exists a conic $C_{x,y}$ through $x$ and
$y$, cutting $\Tx X$ only at $x$. Thus $\pi_x(C_{x,y})=\pi_x(y)\in
W_x$ is a general point and clearly
$\widetilde{\pi}_x(\p(\mathbf{T}_xC_{x,y}))=\pi_x(C_{x,y})=\pi_x(y)$.
Therefore the restriction of $\widetilde{\pi}_x$ to $E$ is dominant
as a map to $W_x\subseteq\p^{N-n-1}$, yielding
$\dim(|II_{x,X}|)=N-n-1$. In particular $N-n-1=\dim(|II_{x,X}|)\leq
\dim(|\O_{\p^{n-1}}(2)|)=\frac{n(n+1)}{2}-1$ and
$N\leq\frac{n(n+3)}{2}$.

Suppose from now on  $\delta\geq 2$. If $y\in X$ is a general point
and if $C_{x,y}$ is a smooth conic through $x$ and $y$ the point
$\p(\mathbf{T}_xC_{x,y})$ is a general point of
$\p((\mathbf{T}_xX)^*)$,  by Theorem \ref{conics} part (6). Consider
the unique quadric hypersurface $Q_{x,y}$ of dimension $\delta\geq
2$ through $x$ and $y$, the irreducible component through $x$ and
$y$ of the entry locus of a general $p\in \langle x,y\rangle $. Then
$C_{x,y}\subset Q_{x,y}$ and $\Tx C_{x,y}\subset \Tx Q_{x,y}$.
Take a line $L_x$ through $x$ and contained in $Q_{x,y}$, which
can be thought of as a point of $Y_x\subset\p((\mathbf{T}_xX)^*)$.
The plane $\langle L_x,\Tx C_{x,y}\rangle $ is contained in $\Tx Q_{x,y}$ so that
it cuts $Q_{x,y}$ at least in another line $L'_x$, clearly different
from $\Tx C_{x,y}$. Thus $\Tx C_{x,y}$ belongs to the pencil
generated by $L_x$ and $L'_x$, which projectivezed in
$\p((\mathbf{T}_xX)^*)$ means that through the general point
$\p(\mathbf{T}_xC_{x,y})\in\p((\mathbf{T}_xX)^*)$ there passes the
secant line $\langle \p(\mathbf{T}_xL_x),\p(\mathbf{T}_xL'_x)\rangle $ to $Y_x$.
Therefore $Y_x\subset\p((\mathbf{T}_xX)^*)$ is non-degenerate and
the join of $Y_x$ with itself equals $\p((\mathbf{T}_xX)^*)$. For an
irreducible $Y_x\subset\p((\mathbf{T}_xX)^*)$ this means exactly
$SY_x=\p((\mathbf{T}_xX)^*)$. The scheme $\Tx Q_{x,y}\cap Q_{x,y}$
is a quadric cone with vertex $x$ and base a smooth quadric
hypersurface of dimension $\delta-2$. The lines in $\Tx Q_{x,y}\cap
Q_{x,y}$ describe a smooth quadric hypersurface of dimension
$\delta-2$, $\widetilde{Q}_{x,y}\subset
Y_x\subset\p((\mathbf{T}_xX)^*)$, whose linear span
$\langle \widetilde{Q}_{x,y}\rangle =\p^{\delta-1}$ passes through
$r=\p(\mathbf{T}_xC_{x,y})$. Since
$\widetilde{\pi}_x:\p((\mathbf{T}_xX)^*)\map W_x\subseteq\p^{N-n-1}$
is given by a linear system of quadrics vanishing on $Y_x$, the
whole $\p^{\delta-1}$ is contracted by $\widetilde{\pi}_x$ to
$\widetilde{\pi}_x(r)$. The closure of the irreducible component of
$\widetilde{\pi}_x^{-1}(\widetilde{\pi}_x(r))$ passing through $r$
has dimension $n-1-\dim(W_x)=\delta-1$ so that it coincides with
$\langle \widetilde{Q}_{x,y}\rangle =\p^{\delta-1}$. This also shows that $Y_x$ is
an irreducible component of the support of the base locus scheme of
$|II_{x,X}|$ and also that, when irreducible,
$Y_x\subset\p((\mathbf{T}_xX)^*)$ is a $QEL$-manifold of type
$\delta-2$. Indeed in this case $r\in SY_x$ is a general point and
every secant or tangent line to the smooth irreducible variety $Y_x$
passing through $r$ is contracted by $\widetilde{\pi}_x$, since the
quadrics in $|II_{x,X}|$ vanish on $Y_x$. Thus every secant line
through $r$ is contained in $\langle \widetilde{Q}_{x,y}\rangle $ and the entry
locus with respect to $r$ is exactly $\widetilde{Q}_{x,y}$. This
concludes the proof of parts (2) and (3).

Suppose $\delta\geq 3$ and let us concentrate on part (4). Item (a)
follows directly from Barth--Larsen Theorems, \cite{BL}, but we provide a direct
proof using the geometry of $LQEL$-varieties. There are lines
through a general point $x\in X$, for example the ones constructed
from the family of entry loci. Reasoning as in Proposition \ref{Yx},
we get
\[TX_{|L}=\O_{\p^1}(2)\oplus\O_{\p^1}(1)^{m(L)}\oplus\O_{\p^1}^{n-m(L)-1}\]
for every line $L$ through $x$. Thus if such a  line comes from a
general entry locus, we get
\[2+m(L)=-K_X\cdot L=\frac{-K_X\cdot C}{2}=\frac{n+\delta}{2},\]
yielding $m(L)=\frac{n+\delta}{2}-2$.

We define $R_x$ to be the locus of points on $X$ which can be joined
to $x$ by a connected chain of lines whose numerical class is
$\frac{1}{2}[C]$, $C\in \mathcal{C}$ a general conic. By
construction we get $R_x=X$, so that
 the Picard number of $X$ is one by
\cite[IV.3.13.3]{Kollar}. Since the variety  $X$ is simply connected
being rationally connected, see Theorem \ref{conics}, we deduce
$\Num(X)=\NS(X)=\Pic(X)\simeq\mathbb{Z}\langle \O(1)\rangle $. Thus $X\subset\p^N$
is a Fano variety, $Y_x\subset\p^{n-1}$ is equidimensional of
dimension $\frac{n+\delta}{2}-2$ and $m(L)=\frac{n+\delta}{2}-2$ for
every line $L$ through $x$. We claim that $Y_x$ is irreducible.

Indeed, if there were  two irreducible components
$Y_x^1,Y_x^2\subset Y_x\subset\p((\mathbf{T}_xX)^*)$, then
$\dim(Y^1_x)+\dim(Y^2_x)=n+\delta-4\geq n-1$, in contrast to
 Proposition
\ref{Yx}.

The fact that $Y_x\subset\p((\mathbf{T}_xX)^*)$ is a $QEL$-manifold
of type $\delta(X)-2$ such that $SY_x=\p((\mathbf{T}_xX)^*)$ follows
from part (3) above. Therefore all the assertions are now proved.
\end{proof}

\begin{Example}\label{Segre}({\it Segre varieties $X=\p^l\times
\p^m\subset\p^{lm+l+m}$, $l\geq 1$, $m\geq 1$, are $QEL$-manifolds
of type $\delta=2$}) By Proposition \ref{K2}, we know that
$X=\p^l\times \p^m\subset\p^{lm+l+m}$ is a $QEL$-manifold, clearly with $\delta\geq
2$, and we
calculate its type, that is we determine  $\delta(X)$.

The locus of lines through a point $x\in X$ is easily described,
being the union of the two linear spaces of the rulings through $x$,
that is  $Y_x=\p^{l-1}\amalg\p^{m-1}\subset\p^{l+m-1}$. Letting
notation be as in Theorem \ref{conics}, we have $C\equiv L_1+L_2$,
where the lines $L_1$ and $L_2$ belongs to different rulings. Then
$n+\delta=-K_X\cdot C=(-K_X\cdot L_1)+(-K_X\cdot
L_2)=(l-1)+2+(m-1)+2=n+2$, so that $\delta(\p^l\times\p^m)=2$ for
every $l,m\geq 1$.
\end{Example}

\begin{Example}\label{Grassmann}({\it Grassmann varieties of lines $\mathbb{G}(1,r)\subset
\p^{\binom{r+1}{2}-1}$ are $QEL$-manifolds of type $\delta=4$}) It
is well known that $Y_x\subset
\p((\mathbf{T}_x\mathbb{G}(1,r))^*)\simeq\p^{2r-3}$ is projectively
equivalent to the Segre variety $\p^{1}\times\p^{r-2}\subset
\p^{2r-3}$. Moreover, $\mathbb{G}(1,r)\subset
\p^{\binom{r+1}{m+1}-1}$ is a $QEL$-manifold, for example by
Proposition \ref{K2}, and we determine its type $\delta$. 
Take $x,y\in
\mathbb{G}(1,r)$ general. They represent two lines
$l_x,l_y\subset\p^r$, $r\geq 3$, which are skew so that
$\langle l_x,l_y\rangle =\p_{x,y}^3\subseteq\p^r$. The Pl\" ucker embedding of
the lines in $\p^3_{x,y}$ is a
$\mathbb{G}(1,3)_{x,y}\subseteq\mathbb{G}(1,r)$ passing through
$x$ and $y$. Therefore $\delta(\mathbb{G}(1,r))\geq 4$.
Thus
$r-1=\dim(Y_x)=-K_X\cdot L-2$, where $L\subset\mathbb{G}(m,r)$ is an
arbitrary line, yielding $-K_X=(r+1)H$, $H$ an hyperplane section.
Finally  $r+1=\frac{2(r-1)+\delta}{2}$  by Theorem \ref{quadrichetangenti},
that is $\delta=4$.
\end{Example}

\begin{Example}\label{E6variety}({\it Spinor variety $S^{10}\subset\p^{15}$ and
$E_6$-variety $X\subset\p^{26}$ as $QEL$-manifolds}) Let us analyze
the 10-dimensional spinor variety $S^{10}\subset\p^{15}$. It is
scheme theoretically defined by 10 quadratic forms defining a map
$\phi:\p^{15}\map\phi(\p^{15})\subset\p^{9}$. The image
$Q=\phi(\p^{15})\subset\p^9$ is a smooth 8-dimensional quadric
hypersurface  and the closure of every fiber is a $\p^7$ cutting $X$
along a smooth quadric hypersurface, see for example \cite{ESB}. In
particular $\delta(X)=6$ and $X\subset\p^{15}$ is a Fano manifold of
index $i(X)=8=n-2$ such that $\Pic(X)=\mathbb{Z}\langle \O_X(1) \rangle $. It is a
so called Mukai variety with $b_2(X)=1$ and by the above description
it is a $QEL$-manifold of type $\delta=6$. For every $x\in X$ the
variety $X^1:=Y_x(X)\subset\p^9$ is a variety of dimension
$\frac{n+\delta(X)}{2}-2=6$, defined by $\codim(X)=5$ quadratic
equations yielding a dominant map $\widetilde{\pi}_x:\p^9\map\p^4$.
The  general fiber of $\widetilde{\pi}_x$ is a linear $\p^5$ cutting
$Y_x\subset\p^9$ along a smooth quadric hypersurface of dimension $4$;
see Theorem \ref{quadrichetangenti}. It is almost clear (and well known) that
$Y_x\simeq \mathbb{G}(1,4)\subset\p^9$ Pl\" ucker embedded. From $Y_x\subset\p^9$ we can
construct the locus of tangent lines and obtain
$X^2:=Y_x(X^1)\subset\p^5$, the Segre $3$-fold
$\p^1\times\p^2\subset\p^5$; see Example \ref{Grassmann}.

We can begin the process with the 16-dimensional variety
$X=E_6\subset\p^{26}$, a Fano manifold of index $i(X)=12$ with
$b_2(X)=1$ and with $\delta(X)=8$. This   is a $QEL$-manifold of type
$\delta=8$, being the center of a $(2,2)$ special Cremona
transformation, see \cite{ESB} and Proposition
\ref{adjunctionCremona}. By applying the above constructions one
obtains $X^1=Y_x(X)=S^{10}\subset\p^{15}$, see \cite[IV]{Zak1}. One
could also apply \cite{Mukai}, since $X^1\subset\p^{15}$
has dimension 10 and type $\delta=6$ so that it is a Fano manifold
of index $i(X)=(n+\delta)/2=8=n-2$. Hence
$X^2=Y_x(X^1)=\mathbb{G}(1,4)\subset\p^9$ and finally
$X^3=Y_x(X^2)=\p^1\times\p^2\subset\p^5$.
\end{Example}

The examples discussed above and the results of Theorem
\ref{quadrichetangenti} suggest to iterate the process, whenever
possible, of attaching  to an $LQEL$-manifold of type $\delta\geq 3$ a
non-degenerate  $QEL$-manifold
$Y_x\subset\p^{n-1}$ of type $\delta-2$ such that $SY_x=\p^{n-1}$. 
If $r\geq 1$ is the largest
integer such that $\delta-2r\geq 1$, and if $X\subset\p^N$ is a
$LQEL$-manifold of type $\delta$, then the process can be iterated
$r$ times, obtaining $QEL$-manifolds of type $\delta-2k\geq 3$ for
every $k=1,\ldots,r-1$.

\begin{Definition} Let $X\subset\p^N$ be an $LQEL$-manifold
of type $\delta\geq 3$. Let 
\[r_X=\sup\{r\in \mathbb{N}: \delta\geq 2r+1\}=[\frac{\delta-1}{2}].\]

 For every $k=1, \ldots, r_X-1$, we define inductively
\[X^k=X^k(z_0,\ldots,z_{k-1})=Y_{z_{k-1}}(X^{k-1}(z_0,\ldots,z_{k-2})),\]
where $z_i\in X^i$, $i=0,\ldots, k-1$, is a general point and where
$X^0=X$.
\end{Definition}

The process  is well defined by Theorem \ref{quadrichetangenti} since for
every $k=1,\ldots, r_X-1$, the variety $X^k$ is a $QEL$-manifold of
type $\delta(X^k)=\delta-2k\geq 3$. The $QEL$-manifold~$X^k$~depends
on the choices of the general points $z_0,\ldots,z_{k-1}$  used to
define it. The type and dimensions of the $X^k$'s are well defined
and we are interested in the determination of these~invariants.

The following result is crucial for the rest of the paper. Its proof
is a direct consequence of  part (4), d) of Theorem
\ref{quadrichetangenti}.

\begin{Theorem}\label{fundamental} Let $X\subset\p^N$ be an $LQEL$-manifold of
type $\delta\geq 3$. Then:
\begin{enumerate}

\item For every $k=1, \ldots, r_X,$ the variety
$X^k\subset\p^{\frac{n+(2^{k-1}-1)\delta}{2^{k-1}}-2k+1}$ is a
$QEL$-manifold of type $\delta(X^k)=\delta-2k$, of dimension
$\dim(X^k)=\frac{n+(2^k-1)\delta}{2^k}-2k,$ such that
$SX^k=\p^{\frac{n+(2^{k-1}-1)\delta}{2^{k-1}}-2k+1};$
in particular,
$\codim(X^{k})=\frac{n-\delta}{2^{k}}+1.$

\item  $2^{r_X}$ divides $n-\delta$, that is  $n\equiv \delta$
{\rm (mod $2^{r_X}$)}.
\end{enumerate}
\end{Theorem}

\begin{Remark}\label{paritysecantrem} Much weaker forms of the
Divisibility Theorem were proposed in \cite[Theorem 0.2]{Ohno} after
long computations with Chern classes.

The hypothesis $\delta\geq 3$ is clearly sharp for the congruence
established in part (2) of Theorem \ref{fundamental}, or for its
weaker form proved in part (4) of Theorem \ref{quadrichetangenti}.
Indeed for the Segre varieties
$X_{l,m}=\p^l\times\p^m\subset\p^{lm+l+m}$, $1\leq l\leq m$, of odd
dimension  $n=l+m$ we have $\delta(X_{l,m})=2$; see Example
\ref{Segre}.

It is worthwhile remarking that the above result is not true for
arbitrary smooth secant defective varieties having $\delta(X)\geq 3$
neither in the weaker form of a parity result. One can consider
smooth non-degenerate complete intersections $X\subset\p^N$ with
$N\leq 2n-2$ and such that $n\not\equiv N-1$ (mod $2$). It is easy
to see that for an arbitrary non-degenerate smooth complete
intersection $X\subset\p^N$ with $N\leq 2n+1$, we have $SX=\p^N$. If
$N\leq 2n-2$, then $\delta(X)=2n+1-N\geq 3$ and $\delta\equiv N-1$
(mod $2$).

Infinite series  of secant defective smooth varieties $X\subset\p^N$
of dimension $n$ with $SX\subsetneq\p^N$, $\delta(X)\geq 3$ and such
that $n\not\equiv \delta(X)$ (mod $2^{r_X}$) can be constructed in
the following~way. Take $Z\subset\p^N$  a smooth $QEL$-manifold of
type $\delta\geq 4$~and dimension $n$ such that $SZ\subsetneq\p^N$.
Consider a $\p^{N+1}$ containing the previous $\p^N$ as a
hyperplane, take $p\in\p^{N+1}\setminus\p^N$ and let
$Y=S(p,Z)\subset\p^{N+1}$ be the cone over $Z$ of vertex $p$. If
$W\subset\p^{N+1}$ is a general hypersurface of degree $d>1$, not
passing through $p$, then $X=W\cap Y\subset\p^{N+1}$ is~a smooth
non-degenerate variety of dimension $n$ such that
$SX=S(p,SZ)\subsetneq\p^{N+1}$. Thus
$\delta(X) =\delta(Z)-1=\delta-1\geq 3$ and $n\not\equiv \delta(X)$ (mod~$2^{r_X}$) 
since $n\equiv \delta$ (mod~$2^{r_X}$). Clearly also
$n\not\equiv \delta(X)$ (mod $2$).

One can take, for example,
$Z_n=\mathbb{G}(1,\frac{n}{2}+1)\subset\p^{\frac{n(n+6)}{8}}$,
$n\geq 8$, which are $QEL$-manifolds of dimension $n\geq 8$ and type
$\delta=4$ such that $SZ\subsetneq\p^{\frac{n(n+6)}{8}}$.
\end{Remark}

\section{Some classification results}\label{classification}

In this section we classify various
important classes of $LQEL$-manifolds.  \cite{IR} 
contains the complete classification of $CC$-manifolds  with $\delta(X)\leq
2$ and hence that of $LQEL$-manifolds of type $\delta=1,2$.

 Let
us recall that a non-degenerate smooth projective variety
$X\subset\p^{\frac{3}{2}n}$ is said to be a {\it Hartshorne variety}
if it is not a complete intersection, where as always $n=\dim(X)$. 
 It is worth remarking that there exist Hartshorne
varieties different from the ones described in item ii) and iii)
of Corollary \ref{deltatranensudue} below, see for example \cite[Proposition 1.9]{EinHart}. This
last fact was kindly pointed out to me by Giorgio Ottaviani.

The first  relevant   application  of the Divisibility Property is the following classification of
$LQEL$-manifolds of  type $\delta>\frac{n}{2}$,
which answers  a problem posed in \cite[0.12.6]{KS}.

\begin{Corollary}\label{deltatranensudue}
Let $X\subset\p^N$ be an $LQEL$-manifold of type
$\delta$ with $\frac{n}{2}<\delta<n$. Then $X\subset\p^N$ is
projectively equivalent to one of the following:
\begin{enumerate}
\item[i)] the Segre $3$-fold $\p^1\times\p^2\subset\p^5$;

\item[ii)] the Pl\" ucker embedding $\mathbb{G}(1,4)\subset\p^9$;

\item[iii)] the $10$-dimensional spinor variety
$S^{10}\subset\p^{15}$;

\item[iv)] a general hyperplane section of
$\mathbb{G}(1,4)\subset\p^9$;
\item[v)] a general hyperplane section of $S^{10}\subset\p^{15}$.
\end{enumerate}
In particular,  $\mathbb{G}(1,4)\subset\p^9$ and
$S^{10}\subset\p^{15}$ are the only $LQEL$-manifolds, modulo projective
equivalence,  which are also
Hartshorne varieties.

\end{Corollary}
\begin{proof} By assumption $\delta>0$. If
$\delta\leq 2$, then  $n=3$ and  $\delta=2=n-1$. Therefore $N=5$ and
$X$ is projectively equivalent to the Segre 3-fold
$\p^1\times\p^2\subset\p^5$ by Proposition \ref{nmenouno} below.

From now on we can assume $\delta\geq 3$ and that $X\subset\p^N$ is 
a Fano manifold with $\Pic(X)=\mathbb{Z}\langle \O(1) \rangle$. By Theorem
\ref{fundamental} there exists an integer $m\geq 1$ such that
$n=\delta+m2^{r_X}$ and since $2\delta>n$ by hypothesis, we have
\begin{equation}\label{stimaforte}
\delta>m2^{r_X}.
\end{equation}

Suppose $\delta=2r_X+2$. From $2r_X+2>m2^{r_X}$ it follows $m=1$ and
$r_X\leq 2$. Hence either  $\delta=4$ and $n=6$ and $X\subset \p^N$ is
a Fano manifold of index $i(X)=(n+\delta)/2=5=n-1$ or $\delta=6$ and
$n=10$ and $X\subset \p^N$ is
a Fano manifold as above and of index $i(X)=(n+\delta)/2=8=n-2$. In the first case by \cite[Theorem 8.11]{Fujita} we get case
ii). In the second case we apply \cite{Mukai}, obtaining
case iii).

Suppose $\delta=2r_X+1$. From \eqref{stimaforte} we get
$2r_X+1>m2^{r_X}$ forcing $m=1$ and $r_X=1,2$. Therefore either
$\delta=3$ and $n=\delta+m2^{r_X}=5$; or $\delta=5$ and $n=\delta+m2^{r_X}=9$.
Reasoning as above, we get cases iv) and v).

To prove the last part let us recall that for a non-degenerate
smooth variety $X\subset\p^{\frac{3}{2}n}$ necessarily
$SX=\p^{\frac{3}{2}n}$; see \cite[V.1.13]{Zak1}. Thus
$\delta(X)=\frac{n}{2}+1>\frac{n}{2}$ and applying  the first part we deduce
that we are either in case ii) or iii) or that $\frac{n}{2}+1=\delta=n$, i.e.\ $n=2$. In the last case $X\subset\p^3$ would be a quadric surface
which is a complete intersection and hence not a Hartshorne variety. This concludes the proof.
\end{proof}

Let us recall that a   smooth non-degenerate irreducible variety
$X\subset\p^{\frac{3}{2}n+2}$ of dimension $n$ such that
$SX\subsetneq\p^{\frac{3}{2}n+2}$ is called {\it a Severi variety},
cf. \cite{Zak3}, \cite{Zak1}.

Another interesting  application of Theorem \ref{fundamental} is the
classification of $LQEL$-manifolds of type $\delta=\frac{n}{2}$. For
such varieties we get  immediately that $n=2, 4, 8$ or $16$ and among them we find
 Severi varieties. 
 Indeed, by \cite[IV.2.1, IV.3.1, IV.2.2]{Zak1}, see also \cite{Russo3},
Severi varieties are  $QEL$-manifolds of type
$\delta=\frac{n}{2}$.
Once we know that $n=2,4, 8$ or 16, it is rather simple to classify Severi varieties,
see also \cite[IV.4]{Zak1} and \cite{Lan2}.
For $n=2,4$ the result is classical and well known while in our approach the $n=8$
case follows from the classification of Mukai manifolds, \cite{Mukai}. The less obvious case
is $n=16$ where we apply a result from \cite{Lan3}, see also \cite{Lan4}, after
an easy reduction via Corollary \ref{deltatranensudue}. What is notable, in our opinion, is
not the fact that this proof is short, easy, natural, immediate and almost self-contained but
the perfect parallel between our argument based on the Divisibility Theorem and
some proofs of Hurwitz Theorem on the dimension of composition algebras over
a field such as the one contained in  \cite[V.5.10]{Lam}, see also \cite[Chap. 10. Sec. 36]{Curtis}.
Surely this connection is well known today, see \cite[pg. 89--91]{Zak1}, but the other
proofs of the classification of Severi varieties did not make this parallel so transparent.
Moreover in Remark \ref{HartYx} below we shall explain an interesting relation between our
approach to the classification of Severi varieties and Hartshorne Conjecture on Complete Intersection.

About this result and the word "generalization" we would like to quote Herman Weyl:
"{\it Before you can generalize, formalize and axiomatize, there must be a mathematical substance}",
\cite{Weyl}. There is no doubt that the mathematical substance in this problem
is entirely due to Fyodor  Zak, who firstly brilliantly solved it in \cite{Zak3}.

\begin{Corollary}\label{nminore16} 
Let $X\subset\p^N$ be an
$LQEL$-manifold of type $\delta=\frac{n}{2}$. Then $n=2,4, 8$ or $16$
and $X\subset\p^N$ is projectively equivalent to one of the
following:

\begin{enumerate}
\item[i)] the cubic scroll $S(1,2)\subset\p^4$;

\item[ii)] the Veronese surface $\nu_2(\p^2)\subset\p^5$ or one of its isomorphic projection in $\p^4$;

\item[iii)] the Segre $4$-fold $\p^1\times\p^3\subset\p^7$;

\item[iv)] a general $4$-dimensional linear section $X\subset\p^7$
of $\mathbb{G}(1,4)\subset\p^9$;

\item[v)] the Segre $4$-fold $\p^2\times\p^2\subset\p^8$ or one of its isomorphic projections in $\p^7$;

\item[vi)] a general $8$-dimensional linear section
$X\subset\p^{13}$ of $S^{10}\subset\p^{15}$;

\item[vii)] the Pl\" ucker embedding
$\mathbb{G}(1,5)\subset\p^{14}$ or one of its isomorphic projection
in $\p^{13}$;

\item[viii)] the $E_6$-variety $X\subset\p^{26}$ or one of its isomorphic projection in $\p^{25}$;

\item[ix)] a $16$-dimensional linearly normal rational variety
$X\subset\p^{25}$, which is a Fano variety of index $12$ with
$SX=\p^{25}$, $\defect(X)=0$ and such that  the base locus of $|II_{x,X}|$, $Z_x\subset\p^{15}$,  is
the union of  a $10$-dimensional spinor variety $S^{10}\subset\p^{15}$
with  $C_pS^{10}\simeq\p^{7}$, $p\in\p^{15}\setminus S^{10}$.
\end{enumerate}
In particular, a Severi variety  $X\subset\p^{\frac{3n}{2}+2}$ is
projectively equivalent to a linearly normal variety as in {\rm ii),
v), vii)} or {\rm viii).}
\end{Corollary}
\begin{proof} By assumption $n$ is even. If $n<6$, then $n=2$ or
$n=4$. If $n=2$, the conclusion is well known, see \cite{Sev} or
Proposition~\ref{nmenouno}. If $n=4$, then $\delta=2=n-2$.
If $H$ is a hyperplane section and if $C\in \mathcal{C}$ is a general conic, then
$(K_X+3H)\cdot C=-n-\delta+2n-2=0$ by part (5) of Theorem
\ref{conics}. Suppose $X\subset\p^N$ is a scroll over a curve, which is rational
by Theorem \ref{conics}. Since for a rational normal scroll either
$SX=\p^N$ or $\dim(SX)=2n+1$,  we get $N=\dim(SX)=2n+1-\delta=7$ so that
$X\subset\p^7$ is a rational normal scroll of degree $4$, which is
the case described in iii). If $X\subset\p^N$ is not a scroll
over a curve, $|K_X+3H|$ is generated by global sections, see \cite[Theorem 1.4]{Ionescu}, and since through two general  points of $X$ there passes
such a conic, we deduce $-K_X=3H$. Thus  $X\subset\p^N$ is a
del Pezzo manifold, getting  cases iii), iv) or v)
by \cite[Theorem 8.11]{Fujita}.

Suppose from now on $n\geq 6$,  $\delta=n/2\geq 3$ and hence
that $X\subset\p^N$ is a Fano manifold with $\Pic(X)=\mathbb{Z}\langle \O(1) \rangle$. By Theorem
\ref{fundamental}, $2^{r_X}$ divides $n-\delta=\frac{n}{2}=\delta$
so that $2^{r_X+1}$ divides $n$ and $\delta=\frac{n}{2}$ is even.
By definition of $r_X$, $\frac{n}{2}=2r_X+2$, so that,  for some
integer $m\geq 1$,  $$m2^{r_X+1}=n=4(r_X+1).$$ Therefore either
$r_X=1$ and  $n=8$, or $r_X=3$ and $n=16$. In the first case
we get that $X\subset\p^N$ is a Fano manifold as above and of index
$i(X)=(n+\delta)/2=6=n-2$ and we are in cases vi) and vii) by \cite{Mukai}.
In the remaining cases $Y_x\subset\p^{15}$ is a 10-dimensional $QEL$-manifold of type
$\delta=6$ so that $Y_x\subset\p^{15}$ is projectively equivalent
to  $S^{10}\subset\p^{15}$ by Corollary \ref{deltatranensudue}. Furthermore  $$N-16=\dim(|II_{X,x}|)+1\leq\operatorname{h}^0(\I_{S^{10}}(2))=10.$$
Thus either $N=26$ and $Y_x\simeq S^{10}\subset\p^{15}$ is the base locus of the second fundamental
form or one easily sees that we are in case ix). If  $|II_{X,x}|\simeq |H^0(\I_{S^{10}}(2))|$, we are in case
viii) by \cite{Lan3}, see also \cite{Lan4}.
\end{proof}
 
The next remark  is somehow hidden in the proof above and, at least to our knowledge,
it seems to have not been noticed before.

\begin{Remark}\label{HartYx} The famous Hartshorne Conjecture on Complete Intesections asserts that a smooth $n$-dimensional variety $X\subset\p^N$ is a
complete intersection if $n>\frac{2N}{3}$, see \cite{Hartshorne}. We recalled above that 
the classification of $LQEL$-manifolds of type $\delta=\frac{n}{2}$  was known classically for $n=2, 4$. Thus suppose $n\geq 6$ and that $X\subset\p^N$ is a $LQEL$-manifold of type $\delta=\frac{n}{2}\geq 3$. In particular $N\geq\dim(SX)=\frac{3n}{2}+1$. Then, by Theorem \ref{quadrichetangenti}, $Y_x\subset\p^{n-1}$
is a non-degenerate $QEL$-manifold of dimension $\frac{3n}{4}-2$ and type $\delta=\frac{n}{2}-2\geq 1$, which is not a complete intersection since 
$$
\ho(\I_{Y_x}(2))\geq\dim(|II_{x,X}|)+1=N-n\geq \frac{n}{2}+1>\frac{n}{4}+1=\codim(Y_x,\p^{n-1}).
$$
Moreover we have that
$$\frac{3n}{4}-2=\dim(Y_x)>\frac{2(n-1)}{3}\mbox{ if and only if } n>16, $$
so that the existence of a $LQEL$-manifold of type $\delta=\frac{n}{2}$ and  of dimension $n> 16$ would have produced a counterexample to Hartshorne Conjecture on Complete Intersections. Obviously  the interesting fact in the above proof is that we can prove directly $n=8, 16$, if $n\geq 6$, without invoking Hartshorne Conjecture.
\end{Remark}

We collect below some classification results used previously in
this section, whose proof is straightforward.

\begin{Proposition}\label{nmenouno} Let $X\subset\p^N$ be
 an $LQEL$-manifold of type  $\delta=n-1\geq 1$. Then $n=2$ or
$n=3$, $N\leq 5$ and $X\subset\p^N$ is projectively equivalent to
one of the following:
\begin{enumerate}
\item $\p^1\times\p^2\subset\p^5$ Segre embedded, or one of its
hyperplane sections;

\item the Veronese surface $\nu_2(\p^2)\subset\p^5$ or one of its isomorphic projections in $\p^4$.
\end{enumerate}
\end{Proposition}
\begin{proof}
 Theorem \ref{quadrichetangenti} part(4) yields
$n-1\!=\!\delta\!\leq\! 2$. Thus either $n\!=\!2$ and $\delta\!=\!1$, or $n\!=\!3$ and~$\delta\!=\!2$.
Suppose first $n=2$ and let $C\subset X$ be an irreducible conic coming from a  general entry locus.
Since $-K_X\cdot C=3$ by Theorem~\ref{conics}, we get $C^2=1$ via
adjunction formula. Moreover, $h^1(\O_X)=0$ by Theorem~\ref{conics},
so that $h^0(\O_X(C))=h^0(\O_{\p^1}(1))+1=3$ and $|C|$ is base point
free. The morphism $\phi=\phi_{|C|}:X\to\p^2$ is dominant and birational since $\deg(\phi)=C^2=1$. Moreover,
$\phi$ sends a
conic $C$ into a line so that  $\phi^{-1}:\p^2\map X\subset\p^N$ is
given by a  sublinear system of $|\O_{\p^2}(2)|$ of dimension at least four and the conclusion
for $n=2$ is now immediate. If $n=3$, we get the conclusion by
passing to a hyperplane section, taking into account Proposition
\ref{deltan}.
\end{proof}
\section{Special Cremona transformations of type
$(2,d)$}\label{Cremonasection}

In this section we  apply the Divisibility Theorem and the classification of conic-connected
manifolds contained in  \cite{IR} to the
classification of special Cremona transformations of type $(2,d)$,
$d\geq 3$. Special Cremona transformations of type $(2,2)$ were
classified in \cite{ESB} as those for which the base locus scheme is
one of the four Severi varieties.

\begin{Definition}(cf.\ \cite{SR,ST1,ESB})
A {\it special
Cremona transformation $\phi:\p^N\map\p^N$ of type $(d_1,d_2)$}
is a birational map given  by forms of degree $d_1\geq 2$, whose base
locus scheme is a smooth irreducible variety and whose inverse is given
by forms of degree~$d_2\geq 2$.
\end{Definition}

\begin{Proposition}\label{adjunctionCremona} Let $X\subset\p^N$ be the center of a special
Cremona transformation $\phi:\p^N\map\p^N$ of type $(2,d)$. Then:
\begin{enumerate}
\item[i)] $X\subset\p^N$ is a $QEL$-manifold of type
$\delta=\frac{n+2}{d}-1$ such that $SX\subset\p^N$ is a hypersurface
of degree $2d-1$;

\item[ii)] $X\subset\p^N$ is projectively normal;

\item[iii)] if $X\subset\p^N$ is a scroll over a curve $C$, then
$N=2n+2$ and $g(C)>0$;

\item[iv)] if $N\leq 2n+1$ and if $H$ is a hyperplane section, then $|K_X+(n-1)H|$ is generated by
global sections, unless $n=2$ and $X\subset\p^5$ is projectively
equivalent to the Veronese surface.

\end{enumerate}
\end{Proposition}
\begin{proof}
By \cite[Proposition 2.5, b)]{ESB}, the variety $X\subset\p^N$ is
non-degenerate.
By \cite[Proposition 2.3]{ESB}, $SX\subset\p^N$ is a hypersurface of
degree $2d-1$ and $X\subset\p^N$ is a $QEL$-manifold of type
\[\delta=2n+1-\dim(SX)=2n+1-N+1=2n-N+2=(2n-N+3)-1=\frac{n+2}{d}-1,\]
where the last equality  follows from  \cite[Lemma 2.4, b)]{ESB},
see also \cite{CK1}. Part i) is~proved.

 Part ii) is an improvement of the main theorem of \cite{BEL},
using the same proof and  taking into account that on
$\pi:\Bl_X(\p^N)\to\p^N$ the linear system
$|\pi^*(\O_{\p^N}(2))-E|,$ $E$ being the exceptional divisor, is
generated by global sections and big, see also \cite{Pensieri}. This
yields $\operatorname{H}^i(\I_X(k))=0$ for every $i>0$ and for every
\[k\geq 2(N-n)-N-1=N-2n-1.\] As $SX\!\subset\!\p^N$ is a
hypersurface, $N\!=\!\dim(SX)\!+\!1=2n\!+\!2\!-\!\delta\!\!\leq 2n\!+\!2$;  part ii)
follows.

Let us consider part iii). Since $SX\subsetneq\p^N$ is a
hypersurface by part i), the secant variety does not fill the
ambient space. For linearly normal scrolls over curves
$\dim(SX)=\min\{N,2n+1\}$, so that  $N=2n+2$. If $g(C)=0$, then
$X\subset\p^{2n+2}$ would be a rational normal scroll for which
$h^0(\I_X(2))=\frac{(n+2)(n+3)}{2}>2n+3$, in contrast with
$h^0(\I_X(2))=2n+3$, \cite[Proposition 2.5]{ESB}.

By \cite[Theorem 1.4]{Ionescu}, the linear system $|K_X+(n-1)H|$ is
spanned by global sections, unless $(X,H)$ is either $(\p^n,\O(1))$,
$(Q^n,\O(1))$, $(\p^2,\O(2))$ or a scroll over a curve. If
$X\subset\p^N$ is the center of a special Cremona transformation,
then it is linearly normal by part i) and it has codimension at
least 2, so that $|K_X+(n-1)H|$ is generated by global sections
unless $X\subset\p^N$ is a scroll over a curve or it is projectively
equivalent to the Veronese surface in $\p^5$. If $N\leq 2n+1$, the
first case is excluded by part iii), while the last case exists
since quadrics through a Veronese surface define a special Cremona
transformation of $\p^5$, see  Proposition \ref{K2}.
\end{proof}

We are in position to prove a general result on  special Cremona
transformations of type $(2,d)$, $d\geq 3$. It  says that centers
$X\subset\p^N$ of such Cremona transformations are very rare,
that  most of them are $QEL$-manifolds of type $\delta=0$ and that there
are no examples with $\delta=1$. For $d\geq
3$ no example of special Cremona transformation of type $(2,d)$
whose center has $\delta\geq 2$ is known, see also Remark \ref{Cremona}.

\begin{Theorem}\label{boundCremona} Let $X\subset\p^N$ be a $QEL$-manifold of type $\delta$
and dimension $n$, which is the center of a special Cremona
transformation $\phi:\p^N\map\p^N$ of type $(2,d)$, $d\geq 3$.

\begin{enumerate}

\item[(i)] If $\delta=\frac{n+2}{d}-1\geq 3$, then $\delta$,  $n$
and $d$ are even numbers.

\item[(ii)] For odd $d\geq 3$,  $N=2d-2$, $n=d-2$
and  the center of $\phi$ is a $QEL$-manifold $X\subset\p^{2d-2}$ of
type $\delta=0$.

\item [(iii)] For even $d\geq 4$, either $N=2d-2$, $n=d-2$
and  the center of $\phi$ is a $QEL$-manifold $X\subset\p^{2d-2}$ of
type $\delta=0$ or $X\subset\p^N$ is a Fano manifold of even
dimension $n$ with $\Pic(X)\simeq\mathbb{Z}\langle \O(1)\rangle $ and index
$\frac{(n-1)(d+1)+3}{2d}$, which is a $QEL$-manifold of even type
$\delta=\frac{n+2}{d}-1\geq 2$.

\end{enumerate}
\end{Theorem}
\begin{proof} Suppose $\delta\geq 3$ and $\delta=2r_X+1$, $r_X\geq
1$. Then from $2r_X+1=\delta=\frac{n+2}{d}-1$, the last equality
coming from Proposition \ref{adjunctionCremona}, we get
$n=2(d(r_X+1)-1)$ so that $n$ is even, contradicting the
Divisibility Theorem.

Therefore $\delta=2r_X+2$ is even as soon as $\delta\geq 3$ and $n$
is also even by the Divisibility  Theorem. From
$2r_X+2=\frac{n+2}{d}-1$, we get $d(2r_X+3)=n+2$ so that $2$ divides
$d$ since it divides $n$ and since $2r_X+3$ is odd. Part i) is
proved.

Suppose $\delta=2$, that is $n=3d-2$. By part i) of Proposition
\ref{adjunctionCremona}, $X\subset\p^N$ is linearly normal, so that
\cite[Theorem 2.2]{IR} yields that  either $X\subset\p^N$ is
projectively equivalent to $\p^l\times\p^{n-l}\subset\p^{(n-l)l+n}$
Segre embedded or $X\subset\p^N$ is a Fano manifold with
$\Pic(X)\simeq\mathbb{Z}\langle \O(1)\rangle $ and index $\frac{n+2}{2}$. The
first case is excluded because  $SX\subset\p^{(n-l)l+n}$ is a
hypersurface only for  $l=2$ and $n=4$, contradicting $n=3d-2\geq
7$. In the second case $n=3d-2$ is even, forcing $d$ even.

Let us suppose $n+2=2d$, that is $\delta=1$. Then
$X\subset\p^{4d-3}$ would be a linearly normal $QEL$-manifold of even
dimension $2d-2\geq 4$, of type $\delta=1$ and such that $SX\subset
\p^{4d-3}$ is a hypersurface, see Proposition
\ref{adjunctionCremona}.  Therefore  a special
Cremona transformation whose center has type $\delta=1$ cannot exist  by 
\cite[Theorem 2.2]{IR}.
\end{proof}

The following classification results  appear to be new.

\begin{Corollary}\label{special23} Let $X\subset\p^N$ be the center of a special
Cremona transformation of type $(2,3)$. Then $N=4$, $X\subset\p^4$
is projectively equivalent to a linearly normal elliptic curve of
degree $5$, and $\phi:\p^4\map\p^4$ is, modulo a projective
transformation, the quadro-cubic special transformation whose
inverse is special and not defined along a smooth elliptic scroll of
degree $5$ and invariant $e=-1$.
\end{Corollary}

\begin{proof} By Theorem \ref{boundCremona}
 we have to consider only the case $N=4$, $\delta=0$ and
$X\subset\p^4$ a smooth linearly normal curve. Let
$\phi:\p^4\map\p^4$ be such a transformation. Take 3 general
quadrics hypersurfaces $Q_1,Q_2, Q_3\in |H^0(\I_X(2))|$. Then
$Q_1\cap Q_2\cap Q_3=X\cup C$ as schemes, where $C=\phi^{-1}(L)$,
$L\subset\p^4$ a general line. Since $\phi$ is of type $(2,3)$,
$\deg(C)=3$ so that $X\subset\p^5$ is an elliptic normal curve of degree $5$. The
linear system $|H^0(\I_X(2))|$ defines a special Cremona
transformation for example by Proposition \ref{K2}. A nice
geometrical description of $\phi$ and of $\phi^{-1}$  is contained
in \cite[Chapter 8, \S 5]{SR}.
\end{proof}

\begin{Corollary}\label{special25} Let $X\subset\p^N$ be the center of a special
Cremona transformation of type $(2,5)$. Then $N=8$ and
$X\subset\p^8$  is a Fano $3$-fold of degree $13$ and sectional genus $8$,
projection from a point of a degree $14$ and sectional genus $8$ Fano
$3$-fold $Y\subset\p^9$.
\end{Corollary}
\begin{proof} By Theorem \ref{boundCremona}
 we have to consider only the case $N=8$, $\delta=0$ and
$X\subset\p^8$ a smooth linearly normal  $3$-fold defined by nine
quadratic equations. In \cite{Pensieri} it is proved that such a
$3$-fold is as in the statement, concluding the proof.
\end{proof}

Let  $\phi:\p^N\map\p^N$ be a special Cremona transformation 
of type $(d_1,d_2)$, having a base locus $X\subset\p^N$ of dimension $n$. By \cite[Proposition 2.3]{ESB} the locus of $d_1$-secant lines to $X\subset\p^N$, let us say $S_{d_1}X\subset\p^N$,
 is an irreducible hypersurface of  degree $d_1\cdot d_2-1$, clearly contained in $SX\subseteq\p^N$. Thus $N-1=\dim(S_{d_1}X)\leq \dim(SX)=2n+1-\delta(X)$. Moreover, by the Trisecant Lemma,  $N=2n+2$ if and only if $d_1=2$ and $X\subset\p^{2n+2}$ is a  $QEL$-manifold
of type $\delta=0$. Now we consider the case $N=2n+1$.

\begin{Corollary}\label{2npiu1}
Let $X\subset\p^{2n+1}$ be the center of a special
Cremona transformation $\phi:\p^{2n+1}\map\p^{2n+1}$ of type $(d_1,d_2)$. Then 
either $d_1 =d_2= 2$,
 $X\subset\p^5$ is projectively equivalent to the Veronese surface
$\nu_2(\p^2)\subset\p^5$ and the inverse transformation
is special and of the same kind; or $d_1 = d_2 = 3$,  $X\subset\p^3$ is a non-hyperelliptic curve of genus
$3$ and degree $6$ and the inverse transformation is special and of
the same kind.
\end{Corollary}
\begin{proof}  By  \cite[Lemma 2.4]{ESB} (see also \cite{CK1})  we
get the following relation
\[2 +n =d_2[(2-d_1)n + 2].\]
If $d_1\geq 3$, then $n=1$ and $d_1 = 3$ and the conclusion is well known; see for example \cite[VIII.4.3]{SR} or
\cite[Proposition 3.1]{Pensieri}. If $d_1=2$, then $n$ is even and $X\subset\p^{2n+1}$ is a projectively
normal $QEL$-manifold of type $\delta=1$ by Proposition \ref{adjunctionCremona}, so that
we can apply \cite[Theorem 2.2]{IR}.  Quadrics through a Veronese surface in $\p^5$
define a special Cremona~transformation whose inverse is of the same kind. Indeed in this case
$\phi:\p^5\map\p^5$ can be interpreted, modulo projective transformations,
as the map associating to a plane conic its dual conic (see also \cite[p.\ 188]{SR},  \cite[Theorem 2.8]{ESB} or apply directly Proposition~\ref{K2}).
\end{proof}

\begin{Remark}\label{Cremona} The argument used in the proof
of Theorem \ref{boundCremona} also yields a bound $n\leq n_0=n_0(d)$
for the dimension of  special Cremona transformations of type
$(2,d)$, $d$ even, which is less sharp  than in the odd case.
Indeed, if $d\geq 4$ is even and if $\delta=2r_X+2$, then the
Divisibility Theorem implies that $2^{r_X}$ divides $d(2r_X+3)-2$,
bounding $r_X$.

In \cite[\S 5]{HKS}   a series of Cremona
transformations is built $\phi_d:\p^{2(d-1)}\map\p^{2(d-1)}$ of type $(2,d)$
for any $d\geq 2$. They are not defined along a scheme $X_d$ of
dimension $d-2$, which is irreducible for $d\geq 3$. For $d=3,4,5$
this scheme is known to be smooth so that for $d=3,5$ there exists
examples of special Cremona transformations as in case (ii) of the
above theorem, as we saw in the previous corollaries. The cases $d=3,4$ are classical and were described
and studied in \cite{SR,ST1,ST2}. Thus the above
theorem cannot be improved and there could exist examples of special
Cremona transformations of type $(2,d)$ with $\delta=0$ also for
$d\geq 6$.

The Severi varieties of dimension $8$ and $16$ yield examples of
special Cremona transformations of type $(2,2)$ with centers of type
$\delta\geq 3$. One could ask if for $d\geq 4$  there exists a
special Cremona transformation of type $(2,d)$ with even $\delta>0$
and whose center is a Fano variety of the kind described in 
Theorem \ref{boundCremona}. The known examples point towards a negative answer to this
question, so that one can conjecture that also for even $d\geq 4$
necessarily $n=d-2$ and $N=2(d-1)$.
 It could be also true that the above
series are the only possible examples for every (odd) $d\geq 3$, the first interesting
cases to be considered being $d=4$ and $d=6$, due to  Corollaries \ref{special23} and \ref{special25}.
\end{Remark}

\section*{ Acknowledgements} 
It is a great pleasure to thank
  Paltin
Ionescu  for interesting discussions on the subject and for many
useful suggestions.

\def\bibaut#1{{\sc #1}}

\end{document}